\documentclass[a4wide,11pt]{scrartcl}

\usepackage[USenglish]{babel}
\usepackage[T1]{fontenc}
\usepackage[utf8]{inputenc}

\usepackage{graphicx,subfigure,tikz,pgfplots,epstopdf}
\graphicspath{{figs/}}

\usepackage{amsmath,amssymb,amsfonts,amsthm,array,geometry,stmaryrd,booktabs,paralist,todonotes,mathtools,pdflscape,siunitx,placeins}

\usepackage{ifpdf}
\ifpdf
\usepackage[pdftex,colorlinks=true,linkcolor=blue,citecolor=green,urlcolor=blue,bookmarks]{hyperref}
\else
\usepackage[dvipdfm,colorlinks=true,linkcolor=blue,citecolor=green,urlcolor=blue,bookmarks]{hyperref}
\fi 

\numberwithin{equation}{section}

\newcommand{\R}{\mathbb{R}}

\newcommand{\bb}[1]{\boldsymbol{#1}}

\begin{document}

\title{The Shizuta - Kawashima Condition for the Barotropic SHTC Two Fluid Model}
\author{Ferdinand Thein\footnotemark[1]}
%
%
\date{}
\maketitle
\begin{abstract}
    Recently the barotropic two fluid model belonging to the class of \emph{symmetric hyperbolic thermodynamically compatible} (SHTC) systems was studied in detail in \cite{Thein2022}.
    There the question was raised whether the dissipative structure introduced by the source terms satisfies the Shizuta - Kawashima condition.
    This well-known condition is a sufficient criterion for the existence of global smooth solutions of the studied system.
    In this work we exploit the dissipative structure of the system under consideration and verify that the Shizuta-Kawashima condition holds.
\end{abstract}

\renewcommand{\thefootnote}{\fnsymbol{footnote}}
\footnotetext[1]{RWTH Aachen University, Templergraben 55, D-52056 Aachen, Germany.\\
\href{mailto:thein@igpm.rwth-aachen.de}{\textit{thein@igpm.rwth-aachen.de}}}
\renewcommand{\thefootnote}{\arabic{footnote}}

\section{Introduction}\label{sec:intro}
In this work we consider a barotropic two fluid model of the \emph{symmetric hyperbolic thermodynamically compatible} (SHTC) class which was recently discussed in \cite{Thein2022}.
This is a submodel of the conservative SHTC model of compressible two-phase flows introduced in \cite{Romenski2009,Romenski2007}.
The obtained results can be useful for a qualitative understanding of the studied model, of the physical processes occurring in two-phase flows and for a comparative analysis of various models.
For further details on the barotropic model and the class of SHTC systems we refer to \cite{Thein2022} and the references therein.

In particular we want to verify that the \emph{Shizuta-Kawashima condition} holds for the system subject of this study. For details on this condition we refer to \cite{Dafermos2016,Ruggeri2021,Shizuta1985}.
The Shizuta-Kawashima condition is a sufficient criterion for the existence of global classical solutions given that the mathematical entropy of the system is convex and the source is semi-dissipative, cf.
\cite{Dafermos2016,Ruggeri2004}. For studies on systems which violate the condition we exemplarily refer to \cite{Bianchini2022,Mascia2010}. A further detailed study which investigates existence for systems when the
condition fails is presented in \cite{Beauchard2011}.
For the study of a necessary, yet not sufficient, weak condition we refer to \cite{Lou2006}. It still remains open to give a necessary \emph{and} sufficient condition.

The rest of the paper is organized as follows: in Section \ref{sec:model} we present the SHTC system under consideration in this paper, study its eigenstructure and discuss the generalized energy.
In Section \ref{sec:mainpart} we first verify the convexity of the mathematical entropy. Next we show that the system is semi-dissipative and finally verify that the Shizuta-Kawashima condition holds.
The paper closes with some concluding remarks and an outlook to future work in Section \ref{sec:conclusion}.
\section{Governing Equations and Characteristic Fields}\label{sec:model}
The PDE system for compressible two-phase flows was discussed in \cite{Romenski2009,Romenski2007}.
In one dimension the system for barotropic flows was thoroughly discussed in \cite{Thein2022}.
Written in terms of the generalized specific energy potential $\Phi = \Phi(\alpha,c,\rho,w)$ it reads
\begin{subequations}\label{sys:gpr_1d_v1}
    \begin{align} 
        \frac{\partial \alpha\rho}{\partial t} + \frac{\partial \alpha\rho u}{\partial x} &= \xi_1,\label{eq:alpha_balance_v1}\\
        \frac{\partial \rho c}{\partial t} + \frac{\partial (\rho c u + \rho \Phi_w)}{\partial x} &= \xi_2,\label{eq:partial_mass_balance_v1}\\
        \frac{\partial\rho}{\partial t} + \frac{\partial\rho u}{\partial x} &= \xi_3,\label{eq:mix_mass_balance_v1}\\
        \frac{\partial\rho u}{\partial t} + \frac{\partial\left(\rho u^2 + p + \rho w\Phi_w\right)}{\partial x} &= \xi_4,\label{eq:mix_mom_balance_v1}\\
        \frac{\partial w}{\partial t} + \frac{\partial(wu + \Phi_c)}{\partial x} &= \xi_5.\label{eq:rel_velocity_balance_v1}
    \end{align} 
\end{subequations}
Here, $\alpha \equiv \alpha_1$ is the volume fraction of the first phase,
which is connected with the volume fraction of the second phase $\alpha_2$ by the saturation law $\alpha_1+\alpha_2 = 1$ and hence we also write $\alpha_2 \equiv 1 - \alpha$.
Throughout this work we assume that both phases are present, i.e.\ $\alpha \in (0,1)$.
The mixture mass density $\rho$ is connected with the phase mass densities $\rho_1,\rho_2$ by the relation $\rho = \alpha_1\rho_1 + \alpha_2\rho_2$.
The phase mass fractions are defined as $c \equiv c_1 = \alpha_1 \rho_1/\rho,\, c_2 = \alpha_2 \rho_2/\rho$ and it is easy to see that $c_1 + c_2 = 1$. Hence we also may use $c_2 \equiv 1 - c$.
The mixture velocity is given by $u = c_1u_1 + c_2u_2$ and $w=u_1 - u_2$ is the relative phase velocity.
The equations describe the balance law for the volume fraction, the balance law for the mass fraction, the conservation of total mass, the total momentum balance law and the balance for the relative velocity.
The phase interaction is present via algebraic source terms which are proportional to thermodynamic forces.
These source terms are given by
\begin{align}
    \xi_1 \equiv \xi_\alpha = -\frac{\rho}{\tau_\alpha}\partial_\alpha\Phi,\quad
    \xi_2 \equiv \xi_c = -\frac{\rho}{\tau_c}\partial_c\Phi,\quad
    \xi_5 \equiv \xi_w = -\frac{\zeta}{\rho}\partial_w\Phi\label{source_terms}
\end{align}
and describe phase pressure relaxation to the common value $\xi_\alpha$, relaxation to equilibrium of the thermodynamic potentials $\xi_c$ and interfacial friction $\xi_w$.
The coefficients $\tau_\alpha,\tau_c$ characterize the rate of pressure and phase relaxation, $\zeta$ is the coefficient of interfacial friction and all coefficients can depend on parameters of state.
Due to mass and momentum conservation throughout this work we assume $\xi_2 = \xi_3 = \xi_4 = 0$. In equilibrium we have $\partial_\alpha\Phi = 0, \partial_c\Phi = 0$ and $\partial_w\Phi = 0$.
For the multi-dimensional system we refer to \cite{Romenski2009,Romenski2007,Thein2022}.
By introducing a stationary equation and a conservation law for the vorticity the multi-dimensional system may also be written in conservative form and again we refer to \cite{Romenski2009,Romenski2007}.
\subsection{Discussion of the mixture equation of state}
We want to specify the derivatives of the generalized energy $\Phi$.
The mixture equation of state (EOS) is defined as the sum of the mass averaged phase equations of state and the kinematic energy of relative motion
\begin{align}
    \Phi(\alpha, c, \rho, w) &= \phi(\alpha,c,\rho) + c(1 - c)\frac{w^2}{2},\label{eq:mix_eos}\\
    \phi(\alpha,c,\rho) &= c_1\phi_1(\rho_1) + c_2\phi_2(\rho_2) = c\phi_1\left(\frac{c\rho}{\alpha}\right) + (1 - c)\phi_2\left(\frac{(1-c)\rho}{1 - \alpha}\right).\label{eq:mix_int_e}
\end{align}
where $\phi_i(\rho_i)$ is the specific energy potential of the $i$-th phase and is assumed to be known.
Note that we have introduced $\Phi,\phi$ and $\phi_i$ for notational convenience to cover the isentropic and the isothermal case at the same time.
To be precise we have with the specific \emph{internal energy} $E,e,e_i$ and the specific \emph{free energy} $F,f,f_i$
\begin{align*}
    \Phi = \begin{cases} E,\,\text{isentropic}\\ F,\,\text{isothermal}\end{cases},\;\phi = \begin{cases} e,\,\text{isentropic}\\ f,\,\text{isothermal}\end{cases}\;
    \text{and}\;\phi_i = \begin{cases} e_i,\,\text{isentropic}\\ f_i,\,\text{isothermal}\end{cases}
\end{align*}
We now introduce the \emph{mixture pressure}
\begin{align}
    p = \alpha_1p_1 + \alpha_2p_2,\label{def:mix_p}
\end{align}
the \emph{specific enthalpy} and the \emph{specific Gibbs energy} of phase $i$
\begin{align}
    h_i(\rho_i) &= e_i(\rho_i) + \frac{p_i(\rho_i)}{\rho_i},\label{def:enthalpy}\\
    g_i(\rho_i) &= e_i(\rho_i) - Ts_i(\rho_i) + \frac{p_i(\rho_i)}{\rho_i}.\label{def:gibbs_energy}
\end{align}
%
%
We further abbreviate
\begin{align*}
	\psi_i(\rho_i) = \phi_i(\rho_i) + \frac{p_i}{\rho_i} = \begin{cases}
		h_i(\rho_i),\, \qquad \text{isentropic}\\
		g_i(\rho_i),\, \qquad \text{isothermal}
	\end{cases}.
\end{align*}
The speed of sound $a_i$ of phase $i$ is given by the following relation
\begin{align}
    a_i^2 = \rho_i\dfrac{\partial \psi_i}{\partial \rho_i} = \begin{dcases}
        \rho_i\left(\dfrac{\partial h_i}{\partial\rho_i}\right)_s,\, \qquad \text{isentropic}, \\
        \rho_i\left(\dfrac{\partial g_i}{\partial\rho_i}\right)_T,\, \qquad \text{isothermal}. 
    \end{dcases}\label{def:sound_speed}
\end{align}
%
For a thermodynamically consistent equation of state the speed of sound is well defined and we do not have to differ between the two cases for the mathematical considerations.
For the mixture EOS we have the following derivatives
\begin{align}
    \frac{\partial \Phi}{\partial\alpha} = \frac{\partial \phi}{\partial\alpha},\quad \frac{\partial \Phi}{\partial c} = \frac{\partial \phi}{\partial c} + (1 - 2c)\frac{w^2}{2},\quad
    \frac{\partial \Phi}{\partial\rho} = \frac{\partial \phi}{\partial\rho}\quad\text{and}\quad
    \frac{\partial \Phi}{\partial w} = c(1 - c)w.\label{eqn:mixE_deriv}
\end{align}
The derivatives for the generalized energy $\phi(\alpha,c,\rho) = c_1\phi_1(\rho_1) + c_2\phi_2(\rho_2)$ of the mixture are 
\begin{equation}
	 \frac{\partial \phi}{\partial\alpha} = \frac{p_2 - p_1}{\rho}, \qquad 
	 \frac{\partial \phi}{\partial c}     = \psi_1(\rho_1) - \psi_2(\rho_2), \qquad 
	 \frac{\partial \phi}{\partial\rho}   = \frac{p}{\rho^2}.
	 \label{eqn.der.isentropic} 
\end{equation}
Further we define the mixture speed of sound $a$ by
\begin{align}
    a^2 := \frac{\partial p}{\partial\rho} = ca_1^2 + (1 - c)a_2^2.\label{mix_a}
\end{align}
Applying the previous results and relations the system can be rewritten in the following form:
\begin{subequations}\label{sys:gpr_1d_v2}
    \begin{align} 
        \frac{\partial\alpha_1\rho}{\partial t} + \frac{\partial\alpha_1\rho u}{\partial x} &= \xi_\alpha,\label{eq:alpha_balance}\\
        \frac{\partial\alpha_1\rho_1}{\partial t} + \frac{\partial\alpha_1\rho_1 u_1}{\partial x} &= \xi_c,\label{eq:partial_mass_balance}\\
        \frac{\partial\rho}{\partial t} + \frac{\partial\rho u}{\partial x} &= 0,\label{eq:mix_mass_balance}\\
        \frac{\partial\rho u}{\partial t}
        + \frac{\partial\left(\alpha_1\rho_1 u_1^2 + \alpha_2\rho_2 u_2^2 + \alpha_1 p_1(\rho_1) + \alpha_2 p_2(\rho_2)\right)}{\partial x} &= 0,\label{eq:mix_mom_balance}\\
        \frac{\partial w}{\partial t} + \dfrac{\partial}{\partial x} \left(\dfrac{1}{2}u_1^2 - \dfrac{1}{2}u_2^2 + \psi_1(\rho_1) - \psi_2(\rho_2)\right) &= \xi_w,\label{eq:rel_velocity_balance}
    \end{align} 
\end{subequations}
For the barotropic case the total energy inequality, which serves as mathematical entropy inequality, reads
\begin{align}
    \sum_{i=1}^2\frac{\partial\alpha_i\rho_i\left(\phi_i + \frac{1}{2}u_i^2\right)}{\partial t}
        + \frac{\partial\alpha_i\rho_i u_i\left(\psi_i + \frac{1}{2}u_i^2\right)}{\partial x}\leq 0, \label{ineq:isoST_energy}
\end{align}
For further comparisons we refer to \cite{Dafermos2016,Ruggeri2021}.

The system (\ref{eq:alpha_balance}) - (\ref{eq:rel_velocity_balance}) can be written in the conservative form where the vector of conserved quantities reads
\[
  \bb{W} = (w_1,w_2,w_3,w_4,w_5)^T \equiv \left(\alpha_1\rho,\alpha_1\rho_1,\rho,\alpha_1\rho_1 u_1 + \alpha_2\rho_2 u_2,u_1-u_2\right)^T.
\]
For the flux written in conserved variables we refer to \cite{Thein2022} for the sake of brevity.
%
%
\subsection{Characteristic Fields}
We want to give the characteristic fields of the barotropic system (\ref{sys:gpr_1d_v2}) in equilibrium, i.e. the sources are zero.
Hence the Jacobian of the system in equilibrium is given by
\begin{align}
    \bb{A}(\bb{W}_{Eq}) = \begin{pmatrix}
        u                                                           & 0                                                             & -\alpha_1u            & \alpha_1  & 0\\[8pt]
        0                                                           & u                                                             & -c_1u                 & c_1       & c_1c_2\rho\\[8pt]
        0                                                           & 0                                                             & 0                     & 1         & 0\\[8pt]
        -a_1^2\dfrac{c_1}{\alpha_1} + a_2^2\dfrac{c_2}{\alpha_2}    & a_1^2 - a_2^2                                                 & -u^2 + c_1a_1^2 + c_2a_2^2\dfrac{\rho}{\rho_2}\dfrac{\partial\rho_2}{\partial w_3}  & 2u       & 0\\[8pt]
        -\dfrac{a_1^2}{\alpha_1\rho} - \dfrac{a_2^2}{\alpha_2\rho}  & \dfrac{a_1^2}{\alpha_1\rho_1} + \dfrac{a_2^2}{\alpha_2\rho_2} & \dfrac{a_1^2}{\rho} - \dfrac{a_2^2}{\rho}\dfrac{\rho}{\rho_2}\dfrac{\partial\rho_2}{\partial w_3}   & 0         & u
    \end{pmatrix}\label{eq:jacobian_prim_var}
\end{align}
Here we have calculated the Jacobian with respect to the conserved quantities and then replaced them with the primitive variables for better readability. Moreover, we abbreviate
\[
  \frac{\partial\rho_2}{\partial w_3} = \frac{w_3(w_3 - w_1) - w_1(w_3 - w_2)}{(w_3 - w_1)^2} \equiv \frac{\alpha_2 - \alpha_1c_2}{\alpha_2^2}.
\]
The eigenvalues in equilibrium can be computed as 
\begin{align}
    \lambda_{1\pm} = u \pm a_1,\quad\lambda_C = u,\quad\lambda_{2\pm} = u \pm a_2\label{sys_eigenvalues}
\end{align}
and we have (up to scaling) the following right eigenvectors
%
\begin{align}
    \begin{split}
        \bb{R}_{1\pm} &= \left(\alpha_1, 1, 1, \lambda_{1\pm}, \pm\frac{a_1}{\alpha_1\rho_1}\right)^T,\quad \bb{R}_{2\pm} = \left(\alpha_1, 0, 1, \lambda_{2\pm}, \mp\frac{a_2}{\alpha_2\rho_2}\right)^T,\\
        \bb{R}_C &= \left(\alpha_1\varepsilon + \dfrac{\alpha_1}{c_1}, 1, \varepsilon, u\varepsilon, 0\right)^T,\quad\varepsilon = -\frac{\alpha_1c_2 - \alpha_2c_1}{\alpha_2c_1}.
    \end{split}\label{sys_eigenvectors_prim}
\end{align}
The subscripts $1,2$ refer to the corresponding phases and these fields are genuine nonlinear.
The pair $(\lambda_C, \bb{R}_C)$ forms a linear degenerate field and hence is associated to contact waves. The detailed calculations can be found in \cite{Thein2022}.
It should be noted that each character of the present fields is independent of the flow.
\section{Convexity, Dissipativity and the Shizuta-Kawashima Condition}\label{sec:mainpart}
\subsection{Convexity}
The convexity of the total energy with respect to the conservative variables is a crucial point to render the system symmetric hyperbolic.
However, due to the complexity no result is available in the literature so far for this system. Here we will give a first attempt under the assumption of equilibrium.
By continuity this gives convexity for states near the equilibrium.
The total energy written in  conservative variables reads
\begin{align*}
    \mathcal{E}(\mathbf{W}) = w_3\Phi(w_1,w_2,w_3,w_5) + \frac{1}{2}\frac{w_4^2}{w_3}.
\end{align*}
Instead of discussing convexity in terms of the conserved quantities we use the results given in \cite{Godrom2003} (\emph{Chap.\ IV, pp.\ $174$}) and study the convexity of
\[
  \mathcal{E}(\alpha,c,v,u,q) = \Phi(\alpha,c,v,q) + u^2/2
\]
with $v$ being the specific volume and $q := vw$.
In particular we have the following mapping $\mathcal{T}:\R^5 \to \R^5$ between these two sets of variables
\begin{align*}
    \bb{W} &= \mathcal{T}(\alpha,c,v,u,q) = \left(\frac{\alpha}{v},\frac{c}{v},\frac{1}{v},\frac{u}{v},\frac{q}{v}\right)^T,\\
    \text{with}\quad (\alpha,c,v,u,q)^T &= \mathcal{T}(\bb{W}) = \left(\frac{w_1}{w_3},\frac{w_2}{w_3},\frac{1}{w_3},\frac{w_4}{w_3},\frac{w_5}{w_3}\right)^T
    = \mathcal{T}(\mathcal{T}(\alpha,c,v,u,q)).
\end{align*}
First we study the generalized energy $\phi(\alpha,c,v)$ and have for the second order derivatives
\begin{align}
    \begin{split}
        \frac{\partial^2 \phi}{\partial\alpha^2} &= v\left(\frac{a_1^2}{\alpha v_1} + \frac{a_2^2}{(1 - \alpha)v_2}\right),\;
        \frac{\partial^2 \phi}{\partial\alpha\partial c} = -\left(\frac{a_1^2}{\alpha} + \frac{a_2^2}{1 - \alpha}\right),\\
        \frac{\partial^2 \phi}{\partial\alpha\partial v} &= -p_1\left(1 - \frac{a_1^2}{v_1p_1}\right) + p_2\left(1 - \frac{a_2^2}{v_2p_2}\right),\\
        \frac{\partial^2 \phi}{\partial c^2} &= \frac{a_1^2}{c} + \frac{a_2^2}{1 - c},\;
        \frac{\partial^2 \phi}{\partial c\partial v} = -\frac{a_1^2 - a_2^2}{v},\\
        \frac{\partial^2 \phi}{\partial v^2} &= \frac{c a_1^2 + (1 - c)a_2^2}{v^2} =: \frac{a^2}{v^2}.
    \end{split}\label{eqn:mix_e_2deriv}
\end{align}
We want to investigate the Hessian of $\phi(\alpha,c,v)$ which is given by
\begin{align}
    &\bb{D}^2e(\alpha,c,v) =\label{internalenergy_hessian}\\
    &\begin{pmatrix}
        %
        v\left(\frac{a_1^2}{\alpha v_1} + \frac{a_2^2}{(1 - \alpha)v_2}\right)
        & -\left(\frac{a_1^2}{\alpha} + \frac{a_2^2}{1 - \alpha}\right)
        & -p_1\left(1 - \frac{a_1^2}{v_1p_1}\right) + p_2\left(1 - \frac{a_2^2}{v_2p_2}\right)\\
        -\left(\frac{a_1^2}{\alpha} + \frac{a_2^2}{1 - \alpha}\right)
        & \frac{a_1^2}{c} + \frac{a_2^2}{1 - c}
        & -\frac{a_1^2 - a_2^2}{v}\\
        -p_1\left(1 - \frac{a_1^2}{v_1p_1}\right) + p_2\left(1 - \frac{a_2^2}{v_2p_2}\right)
        & -\frac{a_1^2 - a_2^2}{v}
        & \frac{a^2}{v^2}
    \end{pmatrix}\notag
\end{align}
We obtain for the first two principal minors
\begin{align*}
    \mathcal{H}_1 &= v\left(\frac{a_1^2}{\alpha v_1} + \frac{a_2^2}{(1 - \alpha)v_2}\right) > 0,\\
    \mathcal{H}_2 &= v\left(\frac{a_1^2}{\alpha v_1} + \frac{a_2^2}{(1 - \alpha)v_2}\right)\left(\frac{a_1^2}{c} + \frac{a_2^2}{1 - c}\right)
    - \left(\frac{a_1^2}{\alpha} + \frac{a_2^2}{1 - \alpha}\right)^2\\
    &= \frac{(1 - \alpha)v_2a_1^2 + \alpha v_1a_2^2}{\alpha(1 - \alpha)v_1v_2}\left(\frac{v_1a_1^2}{\alpha} + \frac{v_2a_2^2}{1 - \alpha}\right)
    - \left(\frac{(1 - \alpha)a_1^2 + \alpha a_2^2}{\alpha(1 - \alpha)}\right)^2\\ &= \frac{a_1^2a_2^2(v_1 - v_2)^2}{v_1v_2\alpha(1 - \alpha)}>0.
\end{align*}
For the third principal minor, i.e. the determinant, we yield
\begin{align*}
    \mathcal{H}_3 &= v\left(\frac{a_1^2}{\alpha v_1} + \frac{a_2^2}{(1 - \alpha)v_2}\right)\left(\frac{a_1^2}{c} + \frac{a_2^2}{1 - c}\right)\frac{a^2}{v^2}
    + 2\left(\frac{a_1^2}{\alpha} + \frac{a_2^2}{1 - \alpha}\right)\frac{a_1^2 - a_2^2}{v}\frac{\partial^2 \phi}{\partial\alpha\partial v}\\
    &- v\left(\frac{a_1^2}{\alpha v_1} + \frac{a_2^2}{(1 - \alpha)v_2}\right)\left(\frac{a_1^2 - a_2^2}{v}\right)^2
    - \frac{a^2}{v^2}\left(\frac{a_1^2}{\alpha} + \frac{a_2^2}{1 - \alpha}\right)^2
    - \left(\frac{a_1^2}{c} + \frac{a_2^2}{1 - c}\right)\left(\frac{\partial^2 \phi}{\partial\alpha\partial v}\right)^2
\end{align*}
Note that at equilibrium we have $p = p_1 = p_2$ which simplifies $\frac{\partial^2 \phi}{\partial\alpha\partial v}$. A lengthy yet straight forward calculation shows $\mathcal{H}_3 = 0$.
Thus $\phi$ is positive semi definite and the energy $\phi$ is convex, although not strictly convex. Now for $\Phi$ we have the second order derivatives as follows
\begin{align}
    \begin{split}
        \frac{\partial^2 \Phi}{\partial\alpha^2} &= \frac{\partial^2 \phi}{\partial\alpha^2},\; \frac{\partial^2 \Phi}{\partial\alpha\partial c} = \frac{\partial^2 \phi}{\partial\alpha\partial c},\;
        \frac{\partial^2 \Phi}{\partial\alpha\partial v} = \frac{\partial^2 \phi}{\partial\alpha\partial v},\;\frac{\partial^2 \Phi}{\partial\alpha\partial q} = 0,\\
        \frac{\partial^2 \Phi}{\partial c^2} &= \frac{\partial^2 \phi}{\partial c^2} - \left(\frac{q}{v}\right)^2,\;
        \frac{\partial^2 \Phi}{\partial v\partial c} = \frac{\partial^2 \phi}{\partial v\partial c} - (1-2c)\frac{q^2}{v^3},\;
        \frac{\partial^2 \Phi}{\partial q\partial c} = (1 - 2c)\frac{q}{v^2},\\
        \frac{\partial^2 \Phi}{\partial v^2} &= \frac{\partial^2 \phi}{\partial v^2} + 3c(1-c)\frac{q^2}{v^4},\frac{\partial^2 \Phi}{\partial q\partial v} = -2c(1-c)\frac{q}{v^3}\\
        \frac{\partial^2 \Phi}{\partial q^2} &= c(1 - c)\frac{1}{v^2}.
    \end{split}\label{eqn:mixE_2deriv}
\end{align}
For equilibrium we further have $u = u_1 = u_2$ and thus $q = 0$. Therefore the Hessian is given by
\begin{align}
    &\left.\bb{D}^2\Phi(\alpha,c,v,q)\right|_{Eq} =\notag\\
    &\begin{pmatrix}
        %
        v\left(\frac{a_1^2}{\alpha v_1} + \frac{a_2^2}{(1 - \alpha)v_2}\right)
        & -\left(\frac{a_1^2}{\alpha} + \frac{a_2^2}{1 - \alpha}\right) & \frac{a_1^2}{v_1p_1} - \frac{a_2^2}{v_2p_2} & 0\\
        -\left(\frac{a_1^2}{\alpha} + \frac{a_2^2}{1 - \alpha}\right) & \frac{a_1^2}{c} + \frac{a_2^2}{1 - c} & -\frac{a_1^2 - a_2^2}{v}& 0\\
        \frac{a_1^2}{v_1p_1} - \frac{a_2^2}{v_2p_2} & -\frac{a_1^2 - a_2^2}{v} & \frac{a^2}{v^2} & 0\\
        0 & 0 & 0 & c(1-c)\frac{1}{v^2}
    \end{pmatrix}\label{gen_energy_hessian}
\end{align}
Thus it is immediate to see that $\Phi$ is also convex with respect to $\alpha,c,v,q$ sufficiently close to the equilibrium.
As noted before this transfers to the convexity of $w_3\Phi(w_1,w_2,w_3,w_5)$ according to \cite{Godrom2003}.
For the total energy we finally conclude
%
%
\begin{align}
    &\left.\bb{D}^2\mathcal{E}(\alpha,c,v,u,q)\right|_{Eq} =
    \begin{pmatrix}
        %
        \bb{D}^2e(\alpha,c,v) & \mathbf{0} & \mathbf{0}\\
        \mathbf{0}^T & 1 & 0\\
        \mathbf{0}^T & 0 & c(1-c)\dfrac{1}{v^2}
    \end{pmatrix}\label{tot_energy_hessian}
\end{align}
and hence convexity is immediately verified near equilibrium. The non-strict convexity is a crucial point which remains to be discussed.
However, similar results have been obtained for a different but related system in \cite{Saleh2020} where a different choice of variables was used.
\subsection{Dissipativity}
In the following we show that the system under consideration is semi-dissipative in the sense of \cite{Dafermos2016}.
Therefore we need to show that
\[
  -\nabla_\mathbf{W}\mathcal{E}\cdot\mathbf{\Xi} \geq \varepsilon|\mathbf{\Xi}|,\;\varepsilon > 0
\]
with $\mathbf{\Xi} = (\xi_1,\xi_2,\xi_3,\xi_4,\xi_5)^T$ being the vector of the source terms.
It is straight forward to obtain
\[
  \nabla_\mathbf{W}\mathcal{E} = \left(\partial_\alpha\Phi, \partial_c\Phi,\Phi + \frac{1}{2}u^2 - (\alpha\partial_\alpha\Phi + c\partial_c\Phi + v\partial_v\Phi + u^2 + q\partial_q\Phi),u^2,\partial_q\Phi\right)^T.
\]
Hence we yield
\begin{align}
    -\nabla_\mathbf{W}\mathcal{E}\cdot\mathbf{\Xi} &= \frac{\rho}{\tau_\alpha}\left(\partial_\alpha\Phi\right)^2 + \frac{\rho}{\tau_c}\left(\partial_c\Phi\right)^2 + \frac{\zeta}{\rho}\partial_w\Phi\partial_q\Phi
    = \frac{\tau_\alpha}{\rho}\mathbf{\Xi}_1^2 + \frac{\tau_c}{\rho}\mathbf{\Xi}_2^2 + \frac{\rho^2}{\zeta}\mathbf{\Xi}_5^2\notag\\
    &\geq \underbrace{\min\left\{\frac{\tau_\alpha}{\rho},\frac{\tau_c}{\rho},\frac{\rho^2}{\zeta}\right\}}_{=:\varepsilon > 0}|\mathbf{\Xi}|^2.
\end{align}
Actually it is one of the modelling principles of the SHTC systems that the dissipative sources are constructed this way, cf.\ \cite{Romenski2009,Romenski2007}.
\subsection{Shizuta-Kawashima-Condition}
The Shizuta-Kawashima condition reads
\[
  \left.(D_\mathbf{W}\mathbf{\Xi}\cdot\mathbf{R}_i)\right|_{Eq}\neq 0
\]
for all eigenvectors $i = 1\pm,2\pm,C$.
Since it is sufficient to find non-zero entries we will focus on the first source term $\xi_\alpha$ for the sake of brevity.
All quantities in the source term may depend on state variables, but due to the special product structure and the fact that we evaluate the condition at equilibrium the calculation simplifies.
Thus we calculate
\begin{align*}
    \left.\nabla_\mathbf{W}\xi_\alpha\right|_{Eq}  &= -\frac{\rho}{\tau_\alpha}\nabla_\mathbf{W}\left(\partial_\alpha\Phi\right)
    = -\frac{1}{\tau_\alpha}\begin{pmatrix}a_1^2\frac{c_1}{\alpha_1^2} + a_2^2\frac{c_2}{\alpha_2^2}\\ -\left(\frac{a^2_1}{\alpha_1} + \frac{a^2_2}{\alpha_2}\right)\\-a_1^2\frac{c_1}{\alpha_1} +
    a_2^2\frac{c_2}{\alpha_2}\frac{\rho}{\rho_2}\frac{\partial\rho_2}{\partial w_3}\\ 0\\ 0\end{pmatrix}.
    %
    %
\end{align*}
%
For the linearly degenerated field $\mathbf{R}_C$ it is obvious that we only have to investigate the first product corresponding to $\xi_\alpha$. After some algebraic manipulations we yield
\begin{align}
    -\tau_\alpha\left.\nabla_\mathbf{W}\xi_\alpha\cdot\mathbf{R}_C\right|_{Eq} = \frac{a_2^2}{\alpha_2} \neq 0. \label{k_cond1}
\end{align}
Now we want to investigate the fields $\mathbf{R}_{i\pm}$. Due to the special structure we will again consider the product giving the first component, i.e.
\begin{align}
    -\tau_\alpha\left.\nabla_\mathbf{W}\xi_\alpha\cdot\mathbf{R}_{i\pm}\right|_{Eq} = (-1)^i\frac{a_i^2}{\alpha_i}\neq 0.\label{k_cond2}
\end{align}
Thus the Shizuta-Kawashima condition holds. Note that it was not necessary to study the components related to $\xi_c$ or $\xi_w$ and thus this result holds with or without these source terms being present in the system.
It is a particular case of interest to study the properties of the system when one of the present phases vanishes, i.e.\ $\alpha \to \{0,1\}$. However, the results \eqref{k_cond1} and \eqref{k_cond2} remain valid in these limit situations.
\section{Conclusion}\label{sec:conclusion}
In the present work we verified that the Shizuta-Kawashima condition holds for the barotropic two fluid system belonging to the SHTC class.
For this purpose we also provided a proof of convexity for the energy at equilibrium which is the first time up to the best knowledge of the author.
Future work includes the generalization to the larger mixture system discussed in \cite{Romenski2009,Romenski2007} and verify convexity in non-equilibrium.
Furthermore, the multi-dimensional case needs to be studied in detail due to the presence of the vorticity in the system.
Thus this work contributes to the overall understanding of the class of SHTC systems for mixtures.
\subsection*{Acknowledgments}
\small{
    %
    F.T. is funded by the DFG SPP 2183 \emph{Eigenschaftsgeregelte Umformprozesse}, project 424334423 and gratefully acknowledges the support by the research training group
    \textit{Energy, Entropy and Dissipative Dynamics (EDDy)} of the DFG - project no. 320021702/GRK2326.
    %
}
\bibliographystyle{abbrv}
\bibliography{hyp22_thein_lit}

\begin{thebibliography}{10}

\bibitem{Beauchard2011}
K.~Beauchard and E.~Zuazua.
\newblock Large time asymptotics for partially dissipative hyperbolic systems.
\newblock {\em Arch. Ration. Mech. Anal.}, 199(1):177--227, 2011.

\bibitem{Bianchini2022}
R.~Bianchini and R.~Natalini.
\newblock Nonresonant bilinear forms for partially dissipative hyperbolic
  systems violating the shizuta--kawashima condition.
\newblock {\em Journal of Evolution Equations}, 22(3):63, Jul 2022.

\bibitem{Dafermos2016}
C.~M. Dafermos.
\newblock {\em Hyperbolic Conservation Laws in Continuum Physics}, volume 325
  of {\em Grundlehren der mathematischen Wissenschaften}.
\newblock Springer Berlin Heidelberg, 4th edition, 2016.

\bibitem{Godrom2003}
S.~K. {Godunov} and E.~Romenski.
\newblock {\em {Elements of Continuum Mechanics and Conservation Laws}}.
\newblock Springer US, 2003.

\bibitem{Lou2006}
J.~Lou and T.~Ruggeri.
\newblock Acceleration waves and weak {S}hizuta-{K}awashima condition.
\newblock {\em Rend. Circ. Mat. Palermo (2) Suppl.}, (78):187--200, 2006.

\bibitem{Mascia2010}
C.~Mascia and R.~Natalini.
\newblock On relaxation hyperbolic systems violating the shizuta--kawashima
  condition.
\newblock {\em Archive for Rational Mechanics and Analysis}, 195(3):729--762,
  Mar 2010.

\bibitem{Romenski2009}
E.~Romenski, D.~Drikakis, and E.~Toro.
\newblock Conservative models and numerical methods for compressible two-phase
  flow.
\newblock {\em Journal of Scientific Computing}, 42(1):68, Jul 2009.

\bibitem{Romenski2007}
E.~Romenski, A.~D. Resnyansky, and E.~F. Toro.
\newblock Conservative hyperbolic formulation for compressible two-phase flow
  with different phase pressures and temperatures.
\newblock {\em Quart. Appl. Math.}, 65(2):259--279, 2007.

\bibitem{Ruggeri2004}
T.~Ruggeri and D.~Serre.
\newblock Stability of constant equilibrium state for dissipative balance laws
  system with a convex entropy.
\newblock {\em Quart. Appl. Math.}, 62(1):163--179, 2004.

\bibitem{Ruggeri2021}
T.~Ruggeri and M.~Sugiyama.
\newblock {\em Classical and relativistic rational extended thermodynamics of
  gases}.
\newblock Springer, Cham, 2021.

\bibitem{Saleh2020}
K.~Saleh and N.~Seguin.
\newblock Some mathematical properties of a barotropic multiphase flow model.
\newblock In {\em Second workshop on compressible multiphase flows: derivation,
  closure laws, thermodynamics}, volume~69 of {\em ESAIM Proc. Surveys}, pages
  70--78. EDP Sci., Les Ulis, 2020.

\bibitem{Shizuta1985}
Y.~Shizuta and S.~Kawashima.
\newblock {Systems of equations of hyperbolic-parabolic type with applications
  to the discrete Boltzmann equation}.
\newblock {\em Hokkaido Mathematical Journal}, 14(2):249 -- 275, 1985.

\bibitem{Thein2022}
F.~Thein, E.~Romenski, and M.~Dumbser.
\newblock Exact and numerical solutions of the {R}iemann problem for a
  conservative model of compressible two-phase flows.
\newblock {\em Journal of Scientific Computing}, 93(3):83, Nov 2022.

\end{thebibliography}
\end{document}